\documentclass[12pt]{article}
\usepackage[english]{babel}

\usepackage{amsmath, amssymb}
\usepackage{graphicx}
\usepackage{aliascnt}

\usepackage{color}
\newcommand\real{{\rm I\! R}}
\newcommand\nat{{\rm I\! N}}
\newcommand\cmplx{\;\hbox{\vrule height6.8pt width0.8pt depth-0.1pt
           \kern-3.6pt {\rm C}}}

\newaliascnt{proposition}{lemma}

\aliascntresetthe{proposition}

\newaliascnt{theorem}{lemma}

\aliascntresetthe{theorem}

\newaliascnt{assumption}{lemma}

\aliascntresetthe{assumption}

\newaliascnt{remark}{lemma}

\aliascntresetthe{remark}

\newaliascnt{definition}{lemma}
\newtheorem{definition}[definition]{Definition}
\aliascntresetthe{definition}

\def\bql#1{\begin{equation}\label{#1}}
\def\bq{\begin{equation}}
\def\eq{\end{equation}}
\def\eref#1{{\rm (\ref{#1})}}

\usepackage[bookmarksopen, colorlinks, linkcolor =blue, urlcolor =red, citecolor =red, pagebackref=true, menucolor =blue]{hyperref}
\addto\extrasenglish{
    
}

\usepackage{subcaption}

\begin{document}

\title{On the numerical solution of a hyperbolic inverse boundary value problem in bounded domains}

\author{Roman Chapko\thanks{Faculty of Applied Mathematics and Informatics, Ivan Franko National University of Lviv, 79000 Lviv, Ukraine} and Leonidas Mindrinos\thanks{Faculty of Mathematics, University of Vienna, Oskar-Morgenstern-Platz 1, 1090 Vienna, Austria}
}

\date{ }
\maketitle

\begin{abstract}
We consider the inverse problem of reconstructing the boundary curve of a cavity embedded in a bounded domain. The problem is formulated in two dimensions for the wave equation. We combine the Laguerre transform with the integral equation method and we reduce the inverse problem to a system of boundary integral equations. We propose an iterative scheme that linearizes the equation using the Fr\'echet derivative of the forward operator. The application of special quadrature rules results to an ill-conditioned linear system which we solve using Tikhonov regularization.  The numerical results show that the proposed method produces accurate and stable reconstructions.

\vspace*{1cm} {\em Keywords:} boundary reconstruction; Laguerre transform; modified single layer potential; non-linear boundary integral equation; quadrature rules; Tikhonov regularization.
\end{abstract}

\section{Introduction}
The inverse problem of reconstructing part of a boundary of an object from overdetermined measurements on the accessible part of the boundary has attracted  much attention in different research areas because of its importance in various applications \cite{Cakoni06, Caorsi03, Massa04, Naik08}. This problem is related to the solution of partial differential equations (PDEs) and because of its non-linearity and ill-posedness is rather complicated in both theoretical and numerical aspects.

Most numerical methods for such kind of problems provide iterative methods with regularization techniques.  However, the use of integral equations for the numerical solution of the boundary reconstruction problem is still possible in various ways.  One possibility is to reduce the boundary value problem directly to a system of non-linear integral equations using Green's theorem \cite{Alves07, Cakoni12, KrRu}. Another approach is to reduce the inverse problem for the PDE to a system of non-linear integral equations and then apply iterative methods \cite{Greece, ChIv, ChIv2, Gin19, Yaman09}. 

In the case of time-dependent inverse problems there exist additional difficulties because of the presence of the independent time variable. Clearly the methods described above can be applied also to non-stationary problems \cite{ChJo, ChKr1, ChKr3, ChKrYo, ChMi}.  Here, there exist different variants for the discretization of the problem with respect to time.

In \cite{ChMi} the authors used the Laguerre transform for the semi-discretization of an inverse boundary problem for a parabolic PDE. This resulted to a sequence of inverse boundary problems for an elliptic PDE. Then, a special potential representation of the solution led to a sequence of non-linear integral equations. In this paper, we extend this approach to an  inverse boundary problem for a hyperbolic PDE.





\begin{description}
\item[Problem formulation] The domain $\Omega$ is doubly connected in $\real^2$ with smooth boundary $\Gamma$ of class $C^2.$ We assume that $\Gamma$ consists of two disjoint curves $\Gamma_1$ and $\Gamma_2$, meaning $\Gamma = \Gamma_1 \cup \Gamma_2,$ with $\Gamma_1 \cap \Gamma_2 = \emptyset,$ such that $\Gamma_1$ is contained in the interior of $\Gamma_2$ (see Fig.\ref{fig1}).  
\begin{figure}
\centering
        \includegraphics[width=.4\linewidth]{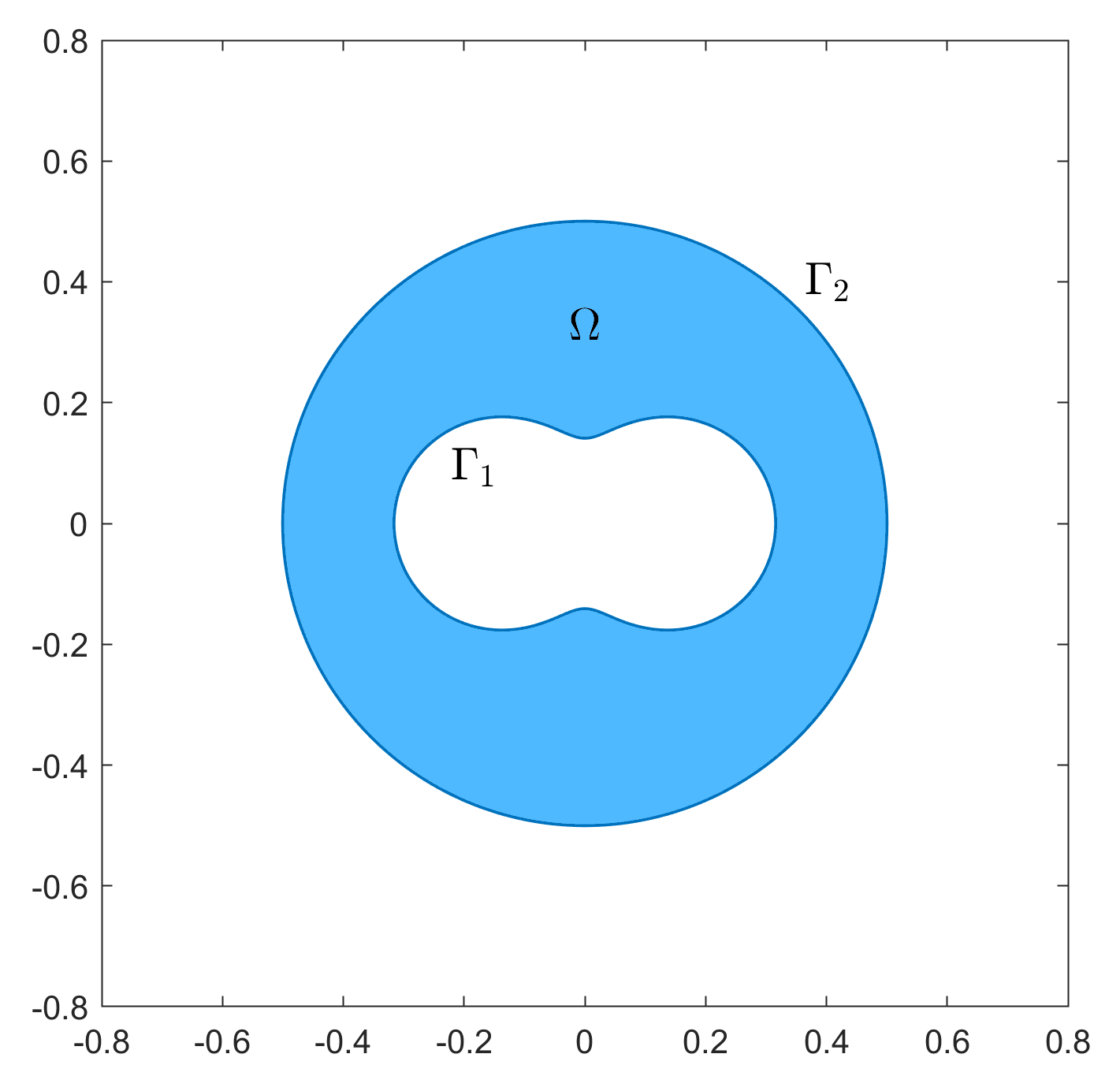}
    \caption{The domain geometry and the notation used throughout this paper.}
\label{fig1}
\end{figure}

We consider the following initial boundary value problem for the wave equation
\bql{1.1}
\frac{1}{a^2}\frac{\partial^2 u}{\partial t^2} =\Delta u, \quad \mbox{in } 
\Omega\times (0,\infty),
\eq
subject to the homogeneous initial conditions 
\bql{1.2}
\frac{\partial u}{\partial t}(\cdot\,,0)=u(\cdot\,,0)=0, \quad\mbox{in } \Omega
\eq
and the boundary conditions
\bql{1.3}
u=0,\quad \mbox{on } \Gamma_1\times[0,\infty), \quad\mbox{and}\quad \frac{\partial u}{\partial \nu}=g \quad \mbox{on }\Gamma_2\times[0,\infty).
\eq
Here $a$ represents the wave speed,  $\nu$ denotes the outward unit normal vector to $\Gamma$ and $g$ is a given and sufficiently smooth function. We refer to \cite{LioMag} for the well-posedness of the direct problem, to find the solution given the domain and the flux $g.$


In this work, we are interested in the numerical solution of the inverse problem to determine the interior boundary curve $\Gamma_1$ from the  knowledge of the Cauchy data on the exterior boundary $\Gamma_2,$ meaning given $g$ and 
\bql{1.4}
u=f, \quad\mbox{on } \Gamma_2 \times (0,\infty).
\eq

\end{description}

An outline of the paper follows: In \autoref{Sec2}, we describe the combination of the Fourier-Laguerre transform with the non-linear boundary integral equation method for the hyperbolic inverse boundary problem.  We derive a sequence of systems of non-linear boundary integral equations, which are transformed into $2\pi$-periodic integral equations. Then, we present an iterative scheme to recover the unknown boundary shape.  

In \autoref{sec3}, we discuss the numerical implementation of the proposed scheme. Given an initial approximation of the unknown boundary curve, we solve the system of equations on the boundary using a quadrature method. The correction of the boundary of the cavity is the solution of the linearized integral equation on the exterior boundary, which we discretize with a trigonometrical collocation method. The Tikhonov regularization is applied to the derived system of linear equations.

Numerical results are presented in \autoref{sec4} confirming that the outlined approach is a feasible way of reconstructing the boundary shape of a cavity.


\section{A two-step approach for dimension reduction}\label{Sec2}
\setcounter{equation}{0}

We first describe the solution $u$ of \eqref{1.1}--\eqref{1.4} using a scaled Fourier expansion
with respect to the Laguerre polynomials. Then, we represent the solution of the stationary problem using a single-layer ansatz. 

\subsection {Semi-discretization in time}
We consider the expansion
$$
 u(x,t)=\kappa\sum_{n=0}^\infty
 u_n(x) {L}_n(\kappa t),
 $$
 where 
$$
  u_n(x)=\int_0^\infty e^{-\kappa t} {L}_n(\kappa t)u(x,t)\,dt,\quad 
n=0,1,2,\ldots
$$
for $\kappa>0,$ using the Laguerre polynomials  ${L}_n$ of order $n$. 

It is easy to show (see for example \cite{ChKr3, ChMi}) that $u$ (sufficiently smooth) is the solution of the time-dependent problem \eqref{1.1}--\eqref{1.4} if and only if its Fourier-Laguerre coefficients $u_n$ satisfy the following sequence of mixed  problems
 \bql{lag10}
 \Delta u_n - \gamma^2 u_n
 =\sum_{m=0}^{n-1} \beta_{n-m}u_m ,
 \quad \mbox{ in }\Omega,
 \eq
 with boundary conditions
 \bql{lag11}
 u_n=0,\quad \mbox{on }\Gamma_1\quad \textrm{and}\quad u_n=f_{n}, \quad\frac{\partial u_n}{\partial \nu}= g_n ,\quad\mbox{on}\;\Gamma_2.
 \eq
 Here $\beta_k=(k+1)\kappa^2/ a^2$, $\gamma^2=\beta_0$ and 
 \begin{align*}
  f_{n}(x) &=\int_0^\infty e^{-\kappa t}
 {L}_n(\kappa t) f(x,t)\,dt,\quad n=0,1,2,\ldots, \\
  g_{n}(x) &=\int_0^\infty e^{-\kappa t}
 {L}_n(\kappa t)g(x,t)\,dt,\quad n=0,1,2,\ldots.
 \end{align*}

In order to apply the non-linear integral equation method we need the sequence of fundamental solutions of the equations \eref{lag10}.

\begin{definition}
  \label{fs-def}
  The sequence of functions $\Phi_n$, for $n=0,1,\ldots$  is called fundamental solution for sequence of equations \eref{lag10}, if it satisfies
  \begin{equation}\label{fs-def-condition}
   \Delta_x \Phi_n(x,y) - \sum_{m=0}^n \beta_{n-m}\Phi_m (x,y)=-\delta(|x - y|). 
  \end{equation}
\end{definition}

We consider the modified Bessel functions
 \bql{besselmod}
 I_0(z)=\sum^\infty_{n=0} \;
 \frac{1}{(n!)^2}\,\left(\frac{z}{2}\right)^{2n},
 \quad
 I_1(z)=\sum^\infty_{n=0} \;
 \frac{1}{n!(n+1)!}\,\left(\frac{z}{2}\right)^{2n+1}
 \eq
 and the modified Hankel functions
\begin{equation}\label{macdonald}
\begin{aligned}
K_0(z)  &= - \left(\ln { \frac{z}2}+C\right)\,I_0(z)
 + \sum^\infty_{n=1}
 \frac{\psi(n)}{ (n!)^2} \,\left(\frac{z}{2}\right)^{2n}, \\
 K_1(z) &= \frac{1}{z}+\left(\ln { \frac{z}2}+C\right)\,I_1(z)
 -\frac{1}{2}\sum^\infty_{n=0}
 \frac{\psi(n+1)+\psi(n)}{ n!(n+1)!}\,\left(\frac{z}{2}\right)^{2n+1}
\end{aligned} 
\end{equation}
  of order zero and one, respectively.  Here, we set $\psi(0)=0$ and
 $$
 \psi(n)=\sum_{m=1}^n\frac{1}{m}\;, \quad n=1,2,\ldots
 $$
 and $C = 0.57721\ldots$ denotes the Euler constant \cite{AbSt}. We define the polynomials $v_n$ and $w_n$ by
 $$
 v_n(r)=
 \sum_{k=0}^{\left[\frac{n}{2}\right]}a_{n,2k}r^{2k},\quad
 w_n(r)=
 \sum_{k=0}^{\left[\frac{n-1}{2}\right]}a_{n,2k+1}r^{2k+1},
 $$
 with the convention $w_0 (r) = 0.$ The coefficients are given by the relations 
\begin{align*}
a_{n,0} &=1, \\
 a_{n,n} &=-\frac{1}{2\gamma n}\;\beta_1 a_{n-1,n-1}, \\
a_{n,k} &=\frac{1}{2\gamma k}
 \left\{4\left[\frac{k+1}{2}\right]^2a_{n,k+1}
 -\sum_{m=k-1}^{n-1}\beta_{n-m} a_{m,k-1}\right\},
 \quad k=n-1,\ldots,1,
\end{align*}  
 for $n=1,2,\ldots$. 
 
Then, following \cite{ChJo,ChKr3}, we  see that the sequence of functions
 \bql{fund-sol}
 \Phi_n(x,y)=K_0(\gamma|x-y|)\,v_n(|x-y|)
 +K_1(\gamma|x-y|)\,w_n(|x-y|),\quad x\neq y,
 \eq
 is a fundamental solution of \eref{lag10} in a sense of \autoref{fs-def}.

\subsection{A boundary integral equation method}

A modified single-layer approach is proposed for solving the sequence of stationary problems. 
We represent the solutions $u_n$ of the problem \eref{lag10} -- \eref{lag11} in the doubly-connected domain $\Omega$  using the following single layer potential form
\bql{pot1}
 u_n(x)=\frac{1}{\pi}\sum_{\ell=1}^2\sum_{m=0}^n
 \int_{\Gamma_\ell} \phi_m^\ell(y) \Phi_{n-m}(x,y)\, ds(y),
 \quad x\in \Omega,
 \eq
with the unknown densities $\phi_m^1$ and $\phi_m^2$, $m=0,1,\ldots$, defined on the boundary curves $\Gamma_1$ and $\Gamma_2$, respectively, and $\Phi_n$ is given by~\eref{fund-sol}.

We let $x$ tend to the boundary $\Gamma$ and using the boundary conditions \eqref{lag11} and the standard jump relations we obtain the following system of equations

\begin{alignat}{3}
\frac1{\pi} \sum_{\ell = 1}^2  \int_{\Gamma_\ell} \phi^\ell_n (y)\Phi_0 (x,y)  ds(y) &= F_{1,n} (x),   &&\quad x &&\in \Gamma_1, \label{eq1}\\
 \phi^2_n (x)+ \frac1{\pi} \sum_{\ell = 1}^2  \int_{\Gamma_\ell}\phi^\ell_n (y) \frac{\partial \Phi_0 }{\partial \nu (x)}  (x,y) ds(y) &= G_n (x),  &&\quad x &&\in \Gamma_2 ,\label{eq2}\\
\frac1{\pi} \sum_{\ell = 1}^2  \int_{\Gamma_\ell} \phi^\ell_n (y)\Phi_0 (x,y)  ds(y) &=  F_{2,n} (x),   &&\quad x &&\in \Gamma_2 ,\label{eq3}
\end{alignat}
for the right-hand sides
\begin{align*}
F_{1,n} (x) &= -  \frac1{\pi} \sum_{\ell = 1}^2 \sum_{m=0}^{n-1} \int_{\Gamma_\ell} \phi^\ell_m (y) \Phi_{n-m} (x,y) ds(y), \\
G_n (x) &= g_n (x)  - \sum_{m=0}^{n-1} \phi^2_m (x)  -  \frac1{\pi} \sum_{\ell = 1}^2 \sum_{m=0}^{n-1} \int_{\Gamma_\ell} \phi^\ell_m (y)\frac{\partial\Phi_{n-m}}{\partial \nu (x)} (x,y)  ds(y), \\
F_{2,n} (x) &= f_n (x) -  \frac1{\pi} \sum_{\ell = 1}^2 \sum_{m=0}^{n-1} \int_{\Gamma_\ell} \phi^\ell_m (y)\Phi_{n-m} (x,y)  ds(y).
\end{align*}

This is a system of three  equations for the three unknowns: the two densities $\phi^1_n , \, \phi^2_n$ and the boundary curve $\Gamma_1.$ The integral operators are singular,  linear on the densities but act non-linearly on the boundary curve. We will consider the Fr\'echet derivative of the integral operators for linearizing them.


Before presenting the iterative method, we consider the parametrization of the system \eqref{eq1}--\eqref{eq3}.  We assume the following parametric representation of the boundary
$$
\Gamma_\ell=\{x_\ell(s)=(x_{1\ell}(s),x_{2\ell}(s)), s\in[0,2\pi]\}, \quad \ell = 1,2
$$  
and we define 
$
\varphi_n^\ell(s) =\phi_n^\ell(x_\ell(s))|x'_\ell(s)|.   
$

Then, the system \eqref{eq1}--\eqref{eq3} takes the form
\begin{alignat}{3}
\frac{1}{2\pi}\sum_{\ell=1}^2\int_0^{2\pi} \varphi_n^\ell(\sigma) H^{1,\ell}_0 (s,\sigma)\, d\sigma  &= \tilde F_{1,n} (s), &&\quad s &&\in [0,2\pi], \label{eq1par}\\
\frac{\varphi_n^2(s)}{|x'_2(s)|}+\frac{1}{2\pi}\sum_{\ell=1}^2\int_0^{2\pi} \varphi_n^\ell(\sigma) Q^{2,\ell}_0 (s,\sigma) d\sigma &= \tilde G_n (s), &&\quad s &&\in [0,2\pi], \label{eq2par}\\
\frac{1}{2\pi}\sum_{\ell=1}^2\int_0^{2\pi} \varphi_n^\ell(\sigma) H^{2,\ell}_0 (s,\sigma)\, d\sigma  &= \tilde F_{2,n} (s), &&\quad s &&\in [0,2\pi], \label{eq3par}
\end{alignat}
for $n=0,\ldots,N$, $N\in\nat,$ with the right-hand sides
\begin{align*}
\tilde F_{1,n} (s) &=  -\frac{1}{2\pi}\sum_{\ell=1}^2\sum_{m=0}^{n-1}\int_0^{2\pi}\varphi_m^\ell(\sigma) H^{1,\ell}_{n-m}(s,\sigma)\, d\sigma, \\
\tilde G_n (s) &= g_{n}(x_2(s))- \frac{1}{|x'_2(s)|}\sum_{m=0}^{n-1} \varphi^2_m (s)
-\frac{1}{2\pi}\sum_{\ell=1}^2\sum_{m=0}^{n-1}\int_0^{2\pi}\varphi_m^\ell(\sigma) Q^{2,\ell}_{n-m}(s,\sigma)\, d\sigma, \\
\tilde F_{2,n} (s) &= f_n (x_2 (s))  -\frac{1}{2\pi}\sum_{\ell=1}^2\sum_{m=0}^{n-1}\int_0^{2\pi}\varphi_m^\ell(\sigma) H^{2,\ell}_{n-m}(s,\sigma)\, d\sigma.
\end{align*}
The kernels are given by
\begin{equation}\label{def_par_ker}
 H^{k,\ell}_ n(s,\sigma)=2\Phi_n (x_k(s),x_\ell(\sigma)), \quad  Q^{k,\ell}_ n(s,\sigma)=2\frac{\partial \Phi_n }{\partial \nu(x)} (x_k(s), x_\ell(\sigma)),
 \end{equation} 
 for $s\neq \sigma, \,\,k,\ell=1,2$,  and $n=0,\ldots,N$. The functions $\Phi_n$ are defined in  \eref{fund-sol}.

\subsection{The iterative scheme}

We solve the derived systems of equations iteratively by splitting them to their well- and ill-posed parts. Following \cite{JohSle}, we first solve the well-posed subsystem to obtain the corresponding densities and then we linearize (with respect to the boundary) the ill-posed subsystem to be solved for the update of the  radial function.

In the following, we assume for simplicity a star-like interior curve with parametrization
\bql{7.4}
x_1(s)=\{r(s) (\cos s,\sin s) : \, s\in  [0,2\pi]\},
\eq
where $r:\real \to(0,\infty)$ is  a  $2\pi-$periodic function representing the radial distance from the origin.

We propose to solve the system of equations \eref{eq1par}--\eref{eq3par} using the iterative scheme: 
\begin{description}
\item[Step 1] Given an initial approximation of $\Gamma_1$, we solve the sequence of well-posed systems of integral equations \eqref{eq1par}--\eqref{eq2par} for $\varphi^1_n, \, \varphi^2_n, \,\, n = 0,...,N.$
\item[Step 2] Keeping now the densities fixed,  we linearize  the ill-posed integral equation \eqref{eq3par} resulting to
\begin{equation}\label{lin_eq}
\sum_{m=0}^{n} \mathcal D_{n-m} [\varphi^1_m, r;q]  (s)  =  f_n (x_2 ( s)) -  \frac1{2\pi} \sum_{\ell=1}^2 \sum_{m=0}^{n}\int_0^{2\pi}  \varphi^\ell_m (\sigma) H^{2,\ell}_{n-m} (s,\sigma)  d\sigma, 
\end{equation}
where $q$ is the radial function of the perturbed boundary. 
We solve the $N$  equations for the radial function $q$ of the perturbed  $\Gamma_1,$ and we update as $r + q.$
\end{description}

The equation \eqref{lin_eq} contains the Fr\'echet derivative $\mathcal D_n$ of the integral operator
with kernel $H^{2,\ell}_n$ with respect to $x_1.$ This is a linear operator on $q$ and its form is obtained by formal differentiation of the kernel $H^{2,\ell}_n$ with respect to $x_1.$ We get  
\[
 \mathcal D_n [\varphi , r;q]  (s) = \frac1{2\pi} \int_0^{2\pi} q (\sigma)   \varphi (\sigma) D_n (s,\sigma) d\sigma,
\]
with kernel
$$
 D_n(s,\sigma) = - \frac{(x_2 (s) -  x_1 (\sigma) ) \cdot (\cos \sigma,\,\sin \sigma)}{|x_2 (s) -  x_1 (\sigma) |}   \tilde{\Phi}_n (|x_2 (s) -  x_1 (\sigma) |),
$$
where
$$
 \tilde{\Phi}_n (r) = K_1 (\gamma r)\,\tilde{v}_n (r)  +K_0 (\gamma r)\,\tilde{w}_n ( r),
$$
for the polynomials
\begin{align*}
\tilde{v}_n(r) &=\gamma \sum_{m=0}^{\left[\frac{n}{2}\right]}a_{n,2m}r^{2m}-2 \sum_{m=1}^{\left[\frac{n-1}{2}\right]}ma_{n,2m+1}r^{2m}, \\
 \tilde{w}_n(r) &=\gamma \sum_{m=0}^{\left[\frac{n-1}{2}\right]}a_{n,2m+1}r^{2m+1}-2 \sum_{m=1}^{\left[\frac{n}{2}\right]}ma_{n,2m}r^{2m-1}.
\end{align*}

Note that the Fr\'echet  derivative operator $\mathcal D_n [\tilde{\varphi} , r;q] $  is injective at the exact solution \cite{ChMi}.

\section{Numerical implementation}\label{sec3}
\setcounter{equation}{0}

The numerical implementation of the iterative scheme has been well examined in \cite{ChMi} for a system similar to \eref{eq1par}--\eref{eq3par}. Thus, in this section we give just a briefly description of it.  We refer to \eqref{eq1par} as the ``field" equations and to \eqref{eq2par} as the ``data" equations.

With the given current approximation of the interior boundary $\Gamma_1$ we consider the ``field" equations \eref{eq1par}. Firstly, we handle the singularity of the parametrized kernels. More precisely, the kernel $H^{\ell,\ell}_n$ in~(\ref{def_par_ker}) admits logarithmic singularity. After lengthy but straightforward calculations, we derive the following decomposition
$$
 H^{\ell,\ell}_n(s,\sigma)= H^{\ell,\ell}_{n,1}(s,\sigma)\ln \left( \frac{4}{e}\sin^2 \frac{s-\sigma}{2}\right)+H^{\ell,\ell}_{n,2}(s,\sigma),
$$
 where
 \begin{eqnarray} 
 H^{\ell,\ell}_{n,1}(s, \sigma) 
& = &-I_0(\gamma|x_\ell (s)-x_\ell (\sigma)|)v_n(|x_\ell (s)-x_\ell (\sigma)|)\nonumber\\[-0.1cm]
\\[-0.1cm]
& & 
+I_1(\gamma|x_\ell (s)-x_\ell (\sigma)|)w_n(|x_\ell (s)-x_\ell (\sigma)|)\nonumber
\end{eqnarray}
 and
$$
 H^{\ell,\ell}_{n,2}(s,\sigma)=H^{\ell,\ell}_n(s,\sigma) - H^{\ell,\ell}_{n,1}(s,\sigma)\ln \left( \frac{4}{e}\sin^2 \frac{s-\sigma}{2}\right)
$$
 with diagonal terms
$$
H^{\ell,\ell}_{n,2}(s, s)=-2C-1-2\ln\left(\frac{\gamma|x'_\ell(s)|}{2}\right)+\frac{2a_{n,1}}{\gamma},\;n=0,1,2,\ldots,N.
$$

Also the kernels $Q_{\ell,\ell}^n$ have logarithmic singularities
$$
 Q^{\ell,\ell}_n(s,\sigma)= Q^{\ell,\ell}_{n,1}(s,\sigma)\ln \left( \frac{4}{e}\sin^2 \frac{s-\sigma}{2}\right)+Q^{\ell,\ell}_{n,2}(s,\sigma),
$$
 where
 \begin{eqnarray} 
 Q^{\ell,\ell}_{n,1}(s, \sigma) & = & 
h^{\ell,\ell}(s,\sigma)\,\{I_1(\gamma|x_\ell(s)-x_k(\sigma)|)
 \tilde v_n(|x_\ell(s)-x_k(\sigma)|)\nonumber\\[-0.1cm]
 \\[-0.1cm]
& & -
 I_0(\gamma|x_\ell(s)-x_k(\sigma)|)
 \tilde w_n(|x_\ell(s)-x_k(\sigma)|)\}\nonumber
\end{eqnarray}
 and
$$
 Q^{\ell,\ell}_{n,2}(s,\sigma)=Q_{\ell,\ell}^n(s,\sigma) - Q^{\ell,\ell}_{n,1}(s,\sigma)\ln \left( \frac{4}{e}\sin^2 \frac{s-\sigma}{2}\right)
$$
 with diagonal terms
$$
 Q^{\ell,\ell}_{n,2}(s, s)=\frac{x^\prime_{\ell,2}(s)x^{\prime\prime}_{\ell,1}(s)-x^\prime_{\ell,1}(s)x^{\prime\prime}_{\ell,2}(s)}
 {|x_\ell^\prime(s)|^3},\quad n=0,1,\ldots,N.
$$
Here we introduced the function 
$$
 h^{\ell,k}(s,\sigma)=\frac{(x_{\ell,1}(s)-x_{k,1} (\sigma))x'_{\ell,2}(s)-(x_{2,\ell}(s)-x_{k,2} (\sigma))x'_{\ell,1}(s)}{|x_k(\sigma)-x_\ell (s)|}.
 $$
Clearly the kernels $H^{k,\ell}_n$ and $Q^{k,\ell}_n$ are smooth for $k\ne\ell$, $k,\ell=1,2$.

Thus, we have to solve the sequence of systems of well-posed $2\pi$ periodical integral equations  \eqref{eq1par} with logarithmic singularities. We use for it the Nystr\"om method with trigonometrical quadrature rules (see for details \cite{ChMi,Kr}).

For the ``data" equations \eqref{eq2par} we apply the collocation method and due to its ill-possedness the received sequence of linear systems is solved by Tikhonov regularization.

\section{Numerical results}\label{sec4}
\setcounter{equation}{0}

We approximate the function $q$ by a trigonometric polynomial of the form 
\begin{equation}\label{eq_radial}
q(s) \approx \sum_{j=0}^{2J} q_j \tau_j (s), \quad \nat \ni J \ll M,
\end{equation}
with
\[
\tau_j (s) =  \left.
  \begin{cases}
    \cos (js), & \text{for } j = 0,...,J ,\\
    \sin ((j-J)s), & \text{for } j = J+1,...,2J.
  \end{cases}
  \right.
\]
We substitute \eqref{eq_radial} in the linearized ``data" equations and at the nodal points $\{s_k\}$ we obtain a linear system, which is ill-posed. We apply Tikhonov regularization. The regularization parameter is chosen initially by trail and error and decreases at every iteration step. 

We simulate the Cauchy data by solving the sequence \eqref{lag10} with boundary conditions
 \[
 u_n= f_{1,n}, \quad \mbox{on }\Gamma_1,\quad \text{and}\quad  u_n = f_{2,n},\quad\mbox{on}\;\Gamma_2,
 \]
for given boundary functions $f_{\ell,n}, \, \ell = 1,2.$ To avoid an inverse crime, we consider double amount of nodal points for the direct problem and afterwards we add noise to the Cauchy data on the boundary $\Gamma_2$ with respect to the $L^2$ norm.  We use the boundary functions
\[
f_{1,n}=0, \quad \text{and} \quad f_{2,n} = \frac{e (2+\kappa n (\kappa (n-1)-4))}{4 (\kappa +1)^{n+3}}, \quad n = 0,...,N.
\]

We consider two examples with different boundary curves:

\begin{description}
\item[Example 1] The interior boundary curve $\Gamma_1$ is a rounded rectangle with radial function
\[
r_1 (s) = (\cos^{10} s  +\sin^{10} s )^{-0.1}
\]
and $\Gamma_2$ is a circle with center $(0,\, 0)$ and radius $1.$

\item[Example 2] Here, both boundary curves are apple-shaped with parametrizations
\[
x_1 (s) = r_1 (s) (\cos s, \, \sin s), \quad \mbox{and} \quad x_2 (s) =  (r_2 (s) \cos s - 0.4, \, r_2 (s) \sin s),
\]
for the radial functions
\[
r_1 (s) = \frac{0.45 + 0.3 \cos s - 0.1 \sin 2s}{1.2 + 0.9 \cos s}, \quad \mbox{and} \quad  
r_2 (s) = \frac{1 + 0.9 \cos s + 0.1 \sin 2s}{0.8 + 0.6 \cos s}.
\]
\end{description}

In both examples, the initial guess is a circle with center $(0,0)$ and radius $r_0.$ We set $\alpha=1$ and $\kappa =1,$ we use $N=10$ Fourier coefficients and we solve at the nodal points with $M=64.$ In the following figures, the brown solid line represents the boundary $\Gamma_2,$ the green dotted line shows the initial guess, the red dashed line is the exact boundary $\Gamma_1$ and its reconstruction is the blue solid line.

\begin{figure}[t]
  \centering
  \includegraphics[width=.9\linewidth]{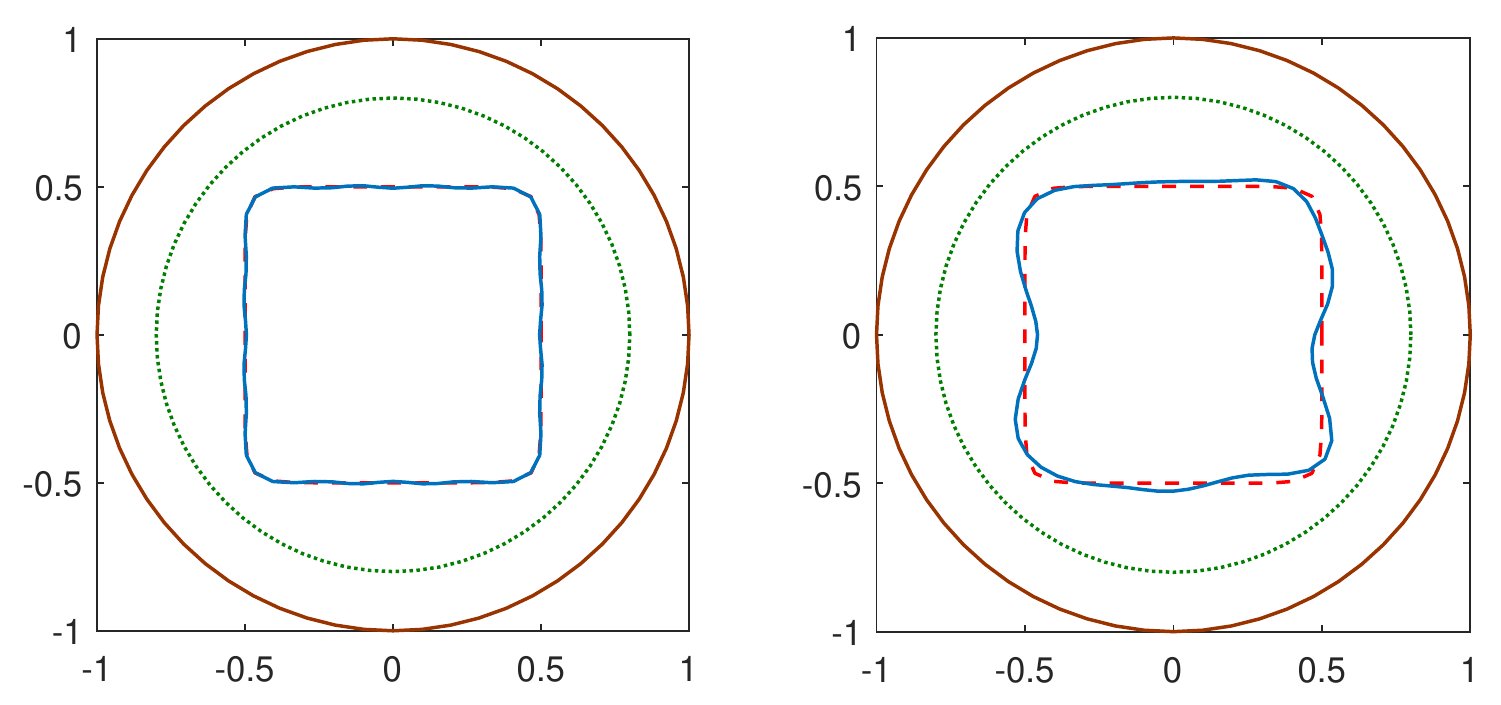}
\caption{Reconstructions of the boundary $\Gamma_1$ of the rounded rectangle for exact data (left) and data with $3\%$ noise (right). }\label{fig_bounded1}
\end{figure}

In the first example, the initial radius is given by $r_0 = 0.8$ and we use $J=13.$ In \autoref{fig_bounded1}, we see the reconstructions for exact (left) and noisy (right) data. The presented results are, with initial regularization parameter $\lambda = 0.01,$ after 21 and 12 iterations, respectively. 

\begin{figure}[h]
  \centering
  \includegraphics[width=.9\linewidth]{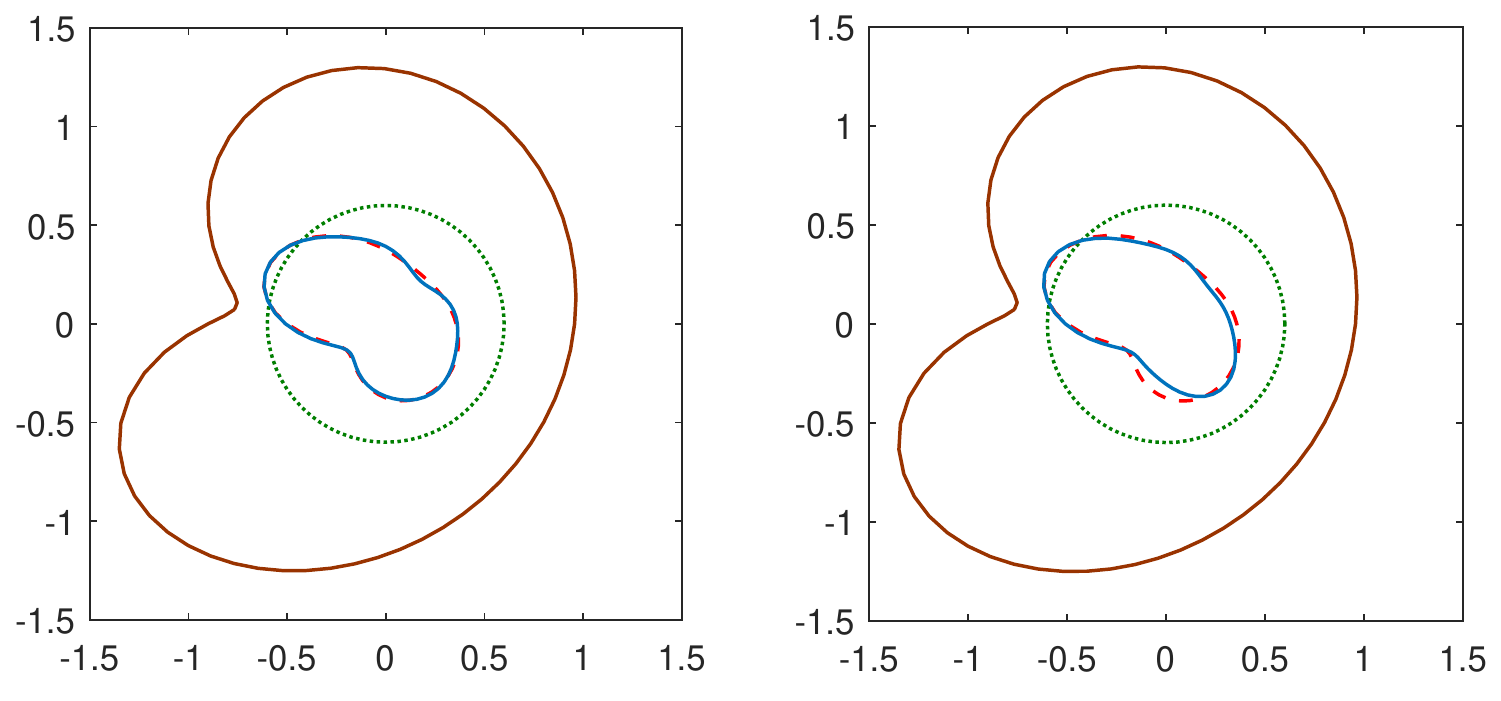}
\caption{Reconstructions of the apple-shaped boundary $\Gamma_1$ for exact data (left) and data with $3\%$ noise (right). }\label{fig_bounded2}
\end{figure}

For the second example, we set $J=5$ and $r_0 = 0.6.$ We consider $\lambda = 0.001$ for the reconstructions presented in \autoref{fig_bounded2}. The algorithm terminated after 10 and 7 iterations, for the noise-free and noisy data, respectively.

We observe that we obtain accurate and relative stable 
reconstructions of the boundary curve. However, we have to stress that the results are sensitive with respect to the initial guess.

\section{Conclusions}
We extended a non-linear integral equations approach for the inverse hyperbolic problem
related to the reconstruction of the interior boundary from the knowledge of
the Cauchy data on the exterior boundary of a doubly connected planar domain. The use of   Laguerre
transform in time leads to a sequence of
stationary inverse boundary problems for the Helmholtz equation.  Next with the help of the
modified single layer potentials these problems were reduced to a sequence of non-linear boundary integral equations. Then, a Newton-type iteration method was applied.  The well-posed system of linear integral equations is solved by
the Nystr\"om method and the ill-posed linear integral equations by the collocation method
with Tikhonov regularization, at every iteration step. Our approach can be extended to the case of three-dimensional domains for similar but more involved fundamental sequences. 

\section{Acknowledgements}
The work of LM was supported by the Austrian Science Fund (FWF) in the project F6801-N36 within the Special Research Programme SFB F68: ``Tomography Across the Scales".

\end{document}